\documentclass[12pt]{article}%
\usepackage{amsmath}
\usepackage{amsfonts}
\usepackage{amssymb}
\usepackage{graphicx}%
\providecommand{\U}[1]{\protect\rule{.1in}{.1in}}
\newtheorem{theorem}{Theorem}

\newtheorem{conjecture}{Conjecture}
\newtheorem{corollary}[theorem]{Corollary}

\newtheorem{problem}{Problem}
\newtheorem{proposition}[theorem]{Proposition}
\newtheorem{question}{Question}
\newtheorem{remark}[theorem]{Remark}

\providecommand{\boksie}{\ensuremath{\mathbin{\raisebox{0.3mm}{$\scriptstyle\square$}}}}
\begin{document}

\title{Dominating and Irredundant Broadcasts in Graphs}
\author{C.M. Mynhardt\thanks{Supported by the Natural Sciences and Engineering
Research Council of Canada.}\\Department of Mathematics and Statistics\\University of Victoria, Victoria, BC, \textsc{Canada}\\{\small kieka@uvic.ca}
\and A. Roux\\Department of Mathematical Sciences\\Stellenbosch University, Stellenbosch, \textsc{South Africa}\\{\small rianaroux@sun.ac.za}}
\date{July 10, 2016}
\maketitle

\begin{abstract}
A \emph{broadcast} on a nontrivial connected graph $G=(V,E)$ is a function
$f:V\rightarrow\{0,1,\dots,\operatorname{diam}(G)\}$ such that $f(v)\leq e(v)$
(the eccentricity of $v$) for all $v\in V$. The cost of $f$ is $\sigma(f)=%
{\textstyle\sum_{v\in V}}
f(v)$. A broadcast $f$ is \emph{dominating} if each $u\in V$ is at distance at
most $f(v)$ from a vertex $v$ with $f(v)\geq1$.

We use properties of minimal dominating broadcasts to define the concept of an
irredundant broadcast on $G$. We determine conditions under which an
irredundant broadcast is maximal irredundant. Denoting the minimum costs of
dominating and maximal irredundant broadcasts by $\gamma_{b}(G)$ and
$\operatorname{ir}_{b}(G)$ respectively, the definitions imply that
$\operatorname{ir}_{b}(G)\leq\gamma_{b}(G)$ for all graphs. We show that
$\gamma_{b}(G)\leq\frac{5}{4}\operatorname{ir}_{b}(G)$ for all graphs $G$.

We also briefly consider the upper broadcast number $\Gamma_{b}(G)$ and upper
irredundant broadcast number $\operatorname{IR}_{b}(G)$, and illustrate that
the ratio $\operatorname{IR}_{b}/\Gamma_{b}$ is unbounded for general graphs.

\end{abstract}

\section{Introduction}

A \emph{broadcast} on a nontrivial connected graph $G=(V,E)$ is a function
$f:V\rightarrow\{0,1,\dots,\operatorname{diam}(G)\}$ such that $f(v)\leq e(v)$
(the eccentricity of $v$) for all $v\in V$. If $G$ is disconnected, we define
a broadcast on $G$ as the union of broadcasts on its components. A broadcast
$f$ is \emph{dominating} if each $u\in V$ is at distance at most $f(v)$ from a
vertex $v$ with $f(v)\geq1$. The \emph{cost} of a broadcast $f$ is
$\sigma(f)=\sum_{v\in V}f(v)$, and the \emph{broadcast number} of $G$ is
\[
\gamma_{b}(G)=\min\left\{  \sigma(f):f\text{ is a dominating broadcast of
}G\right\}  .
\]

A dominating broadcast $f$ such that $f(v)\in\{0,1\}$ for each $v\in V$
corresponds to a \emph{dominating set }of $G$. A dominating set $D$ is
\emph{minimal dominating} (i.e., no subset of $D$ is dominating) if and only
if each $v\in D$ dominates a vertex that is not dominated by $D-\{v\}$, that
is, if and only if $D$ is \emph{irredundant}. Cockayne, Hedetniemi and Miller
\cite{CHM} introduced the concept of irredundance as precisely the property
that makes a dominating set minimal dominating.

Ahmadi, Fricke, Schroeder, Hedetniemi and Laskar \cite{Ahmadi} use a property
that makes a dominating broadcast minimal dominating, which was first
mentioned in \cite{Ethesis}, to define broadcast irredundance, which we state
here in Section \ref{Subsec_Irr}. The \emph{broadcast irredundance number} of
$G$ is defined as%
\[
\operatorname{ir}_{b}(G)=\min\{\sigma(f):f\text{ is a maximal irredundant
broadcast of }G\}.
\]
The definitions imply that $\operatorname{ir}_{b}(G)\leq\gamma_{b}(G)$ for all
graphs $G$, and as our main result we prove that the ratio $\gamma
_{b}/\operatorname{ir}_{b}$ is bounded:

\begin{theorem}
\label{Thm_irb-gammab_bound}For any graph $G$, $\gamma_{b}(G)\leq\frac{5}%
{4}\operatorname{ir}_{b}(G)$.
\end{theorem}

After defining our basic concepts in Section \ref{SecDefs}, we present some
properties of irredundant broadcasts in Section \ref{Sec_Ir}, the most
important of which is a necessary and sufficient condition for an irredundant
broadcast to be maximal irredundant (Theorem \ref{ThmMax_irr}). Theorem
\ref{Thm_irb-gammab_bound} is proved in Section \ref{Sec_Bound}. We briefly
discuss upper broadcast domination and irredundance in Section \ref{SecUpper},
illustrating that the ratio $\operatorname{IR}_{b}/\Gamma_{b}$ is unbounded
for general graphs, and conclude with a list of open problems and conjectures
in Section~\ref{SecOpen}.

\section{Definitions}

This section contains more definitions concerning dominating broadcasts,
neighbourhoods and boundaries of broadcasting vertices, minimal dominating
broadcasts and, finally, irredundant broadcasts. For undefined concepts we
refer the reader to \cite{CLZ}.\label{SecDefs}

\subsection{Dominating Broadcasts}

\label{Sec_Dom}Consider a broadcast $f$ on a connected graph $G=(V,E)$. Define
$V_{f}^{+}=\{v\in V:f(v)>0\}$ and partition $V_{f}^{+}$ into the two sets
$V_{f}^{1}=\{v\in V:f(v)=1\}$ and $V_{f}^{++}=V_{f}^{+}-V_{f}^{1}$. A vertex
$u$ \emph{hears} the broadcast $f$ from some vertex $v\in V_{f}^{+}$, and $v$
$f$-\emph{dominates} $u$, if the distance $d(u,v)\leq f(v)$. Denote the set of
all vertices that do not hear $f$ by $U_{f}$; thus $f$ is a dominating
broadcast if $U_{f}=\varnothing$. A dominating broadcast $f$ of $G$ such that
$\sigma(f)=\gamma_{b}(G)$ is called a $\gamma_{b}$-\emph{broadcast}. Broadcast
domination was introduced by Erwin \cite{Ethesis}, who also gave the trivial
upper bound%
\[
\gamma_{b}(G)\leq\{\gamma(G),\operatorname{rad}(G)\}
\]
for any graph $G$. Graphs for which $\gamma_{b}(G)=\operatorname{rad}(G)$ are
called \emph{radial graphs}. Radial trees are characterized in \cite{Herke,
HM}, where a formula for the broadcast number $\gamma_{b}(T)$ of a tree $T$,
as well as a simple algorithm to determine $\gamma_{b}(T)$, can also be found.
Another algorithm to determine $\gamma_{b}(T)$ is given in \cite{DDH}.

If $f$ and $g$ are broadcasts on $G$ such that $g(v)\leq f(v)$ for each $v\in
V$, we write $g\leq f$. If in addition $g(v)<f(v)$ for at least one $v\in V$,
we write $g<f$. Also, $g\geq f$ ($g>f$, respectively) if $f\leq g$ ($f<g$,
respectively). A dominating broadcast $f$ on $G$ is a \emph{minimal dominating
broadcast} if no broadcast $g<f$ is dominating. Clearly, a $\gamma_{b}%
$-broadcast is a minimal dominating broadcast, but the converse need not be
true. The \emph{upper broadcast number }of $G$, first defined in
\cite{Ethesis} and also studied in \cite{Ahmadi, Dunbar}, is
\[
\Gamma_{b}(G)=\max\left\{  \sigma(f):f\text{ is a minimal dominating broadcast
of }G\right\}  ,
\]
and a dominating broadcast $f$ of $G$ such that $\sigma(f)=\Gamma_{b}(G)$ is
called a $\Gamma_{b}$-\emph{broadcast}. If $f$ is a dominating broadcast such
that $f(v)\in\{0,1\}$ for each $v\in V$, then $\{v\in V:f(v)=1\}$ is a
dominating set of $G$; the smallest cardinality of a dominating set is the
\emph{domination number }$\gamma(G)$, and the largest cardinality of a minimal
dominating set is the \emph{upper domination number }$\Gamma(G)$. Again, Erwin
\cite{Ethesis} gave the trivial lower bound%
\[
\Gamma_{b}(G)\geq\max\{\Gamma(G),\operatorname{diam}(G)\}
\]
for all graphs $G$.

Broadcast domination can be considered as an integer programming (IP) problem.
Its fractional relaxation linear program (LP) has a dual linear program (DLP)
whose IP formulation provides a lower bound for the broadcast number via the
strong duality theorem of linear programming. The dual to the broadcast
domination problem was referred to in \cite{DDH} and studied explicitly by
Teshima \cite{LauraT} and Brewster, Mynhardt and Teshima \cite{BMT}, who
called it the \emph{multipacking problem}. For a positive integer $s$, the
$s$\emph{-neighbourhood} $N_{s}[v]$ of $v\in V$ is the set of all vertices
within distance $s$ from $v$. A set $M$ of vertices of $G$ is called a
\emph{multipacking} if, for each $v\in V$ and each integer $s$ such that
$1\leq s\leq e(v)$, the set $N_{s}[v]$ contains at most $s$ vertices from $M$.
The \emph{multipacking number }$\operatorname{mp}(G)$ is the maximum
cardinality of a multipacking of $G$. The duality of multipackings and
broadcasts implies that $\gamma_{b}(G)\geq\operatorname{mp}(G)$ for any graph
$G$. Hence the existence of a multipacking of cardinality $m$ in a graph $G$
with a dominating broadcast of cost $m$ serves as a certificate that
$\gamma_{b}(G)=m$.

\subsection{Neighbourhoods and boundaries}

For a set $S$ of vertices of a graph $G$ and $s\in S$, the \emph{private
neighbourhood} $\operatorname{PN}(s,S)$ of $s$ with respect to $S$ is the set
of all vertices in the closed neighbourhood of $s$ that are not in the closed
neighbourhood of any other vertex in $S$. If $u\in\operatorname{PN}(s,S)-S$,
then $u$ is an \emph{external private neighbour }of $s$. If $u\in
\operatorname{PN}(s,S)\cap S$, then $u=s$, $s$ is isolated in $G[S]$ and $s$
is said to be an $S$-\emph{self-private neighbour}. For a broadcast $f$ on $G$
and $v\in V_{f}^{+}$, define the

\begin{itemize}
\item $f$-\emph{neighbourhood} of $v$ as $N_{f}[v]=\{u\in V(G):d(u,v)\leq
f(v)\}=N_{f(v)}[v]$

\item $f$-\emph{boundary} of $v$ as\emph{ }$B_{f}(v)=\{u\in
V(G):d(u,v)=f(v)\}$

\item $f$-\emph{private neighbourhood }of\emph{ }$v$ as $\operatorname{PN}%
_{f}(v)=\{u\in N_{f}[v]:u\notin N_{f}[w]$ for all $w\in V^{+}-\{v\}\}$

\item $f$-\emph{private boundary} of $v$ as $\operatorname{PB}_{f}(v)=\{u\in
N_{f}[v]:u$ is not dominated by $(f-\{(v,f(v)\})\cup\{(v,f(v)-1)\}$

\item \emph{external }$f$-\emph{private boundary} of $v$ as
$\operatorname{EPB}_{f}(v)=\operatorname{PB}_{f}(v)-\{v\}$.
\end{itemize}

Note that if $v\in V_{f}^{1}$, then $\operatorname{PB}_{f}%
(v)=\operatorname{PN}_{f}(v)$, and if $v\in V_{f}^{++}$, then
$\operatorname{PB}_{f}(v)=B_{f}(v)\cap\operatorname{PN}_{f}(v)$. Hence the
external $f$-private boundary of $v$ differs from its $f$-private boundary
only if $v\in V_{f}^{1}$ and $v\in\operatorname{PN}_{f}(v)$. As proved in
\cite{BC}, every graph without isolated vertices has a minimum dominating set
in which every vertex has an external private neighbour. A similar result
holds for minimum cost dominating broadcasts.

\begin{proposition}
Every graph without isolated vertices has a $\gamma_{b}$-broadcast $f$ such
that each vertex in $V_{f}^{+}$ has a nonempty external $f$-private boundary.
\end{proposition}

\noindent\textbf{Proof.\hspace{0.1in}}Consider a $\gamma_{b}$-broadcast $f$
such that $|V_{f}^{+}|$ is minimum. Suppose $\operatorname{EPB}_{f}%
(v)=\varnothing$ for some $v\in V_{f}^{+}$. Then every vertex in $B_{f}(v)$
hears a broadcast from some vertex $u\in V_{f}^{+}$. If $v\in V_{f}^{++}$,
then $g=(f-\{(v,f(v)\})\cup\{(v,f(v)-1)\}$ is a dominating broadcast such that
$\sigma(g)<\sigma(f)=\gamma_{b}(G)$, which is impossible. Hence assume
$f(v)=1$. Let $b\in B_{f}(v)=N(v)$ and let $u$ be a vertex in $V_{f}^{+}$ such
that $b$ hears the broadcast from $u$. Then $h=(f-\{(v,f(v)\})\cup
\{(u,f(u)+1)\}$ is a dominating broadcast such that $|V_{h}^{+}|=|V_{f}%
^{+}|-1$, contrary to the choice of $f$.~$\blacksquare$

\subsection{Minimal dominating broadcasts}

The property that makes a dominating broadcast a minimal dominating broadcast,
determined in \cite{Ethesis}\emph{,} is important in the study of the upper
broadcast number $\Gamma_{b}$ and essential for broadcast irredundance. We
restate it here in terms of private boundaries.

\begin{proposition}
\label{PropMinimal}\emph{\cite{Ethesis}}\hspace{0.1in}A dominating broadcast
$f$ is a minimal dominating broadcast if and only if $\operatorname{PB}%
_{f}(v)\neq\varnothing$ for each $v\in V_{f}^{+}$.
\end{proposition}

\noindent\textbf{Proof.\hspace{0.1in}}Suppose $f$ is a minimal dominating
broadcast. For any $v\in V_{f}^{+}$, the broadcast $g=(f-\{(v,f(v)\})\cup
\{(v,f(v)-1)\}$ is not dominating. Hence there either exists a vertex $u\in
B_{f}(v)$ that does not hear the broadcast $g$, in which case $u\in
B_{f}(v)\cap\operatorname{PN}_{f}(v)$, or $f(v)=1$ and $v$ does not hear the
broadcast $g$, in which case $v\in\operatorname{PN}_{f}(v)$. In either case
$\operatorname{PB}_{f}(v)\neq\varnothing$.

Conversely, suppose $\operatorname{PB}_{f}(v)\neq\varnothing$ for each $v\in
V_{f}^{+}$. Then either $v\in V_{f}^{1}\cap\operatorname{PN}_{f}(v)$ and
$g=(f-\{(v,f(v)\})\cup\{(v,f(v)-1)\}$ does not dominate $v$, or $g$ does not
dominate some $u\in\operatorname{PB}_{f}(v)-\{v\}$. Hence $f$ is a minimal
dominating broadcast.~$\blacksquare$

\begin{remark}
A minimal dominating set is a minimal dominating broadcast. Hence $\Gamma
_{b}(G)\geq\Gamma(G)$ for any graph $G$.
\end{remark}

\subsection{Irredundant broadcasts}

\label{Subsec_Irr}Just like irredundance was originally defined in \cite{CHM}
to be precisely the property that makes a dominating set minimal dominating,
Ahmadi et al.~\cite{Ahmadi} define \label{Def_ir}a broadcast $f$ to be
\emph{irredundant} if $\operatorname{PB}_{f}(v)\neq\varnothing$ for each $v\in
V_{f}^{+}$. An irredundant broadcast $f$ is \emph{maximal irredundant} if no
broadcast $g>f$ is irredundant. The \emph{lower} and \emph{upper broadcast
irredundant numbers} of $G$ are
\[
\operatorname{ir}_{b}(G)=\min\left\{  \sigma(f):f\text{ is a maximal
irredundant broadcast of }G\right\}
\]
and
\[
\operatorname{IR}_{b}(G)=\max\left\{  \sigma(f):f\text{ is an irredundant
broadcast of }G\right\}  ,
\]
respectively. Proposition \ref{PropMinimal} and the above definitions imply
the following two results.

\begin{corollary}
\label{Cor_ir-dom}\emph{\cite{Ahmadi}}\hspace{0.1in}Any minimal dominating
broadcast is maximal irredundant.
\end{corollary}

\begin{corollary}
\label{Cor_gammab-irb}\emph{\cite{Ahmadi}}\hspace{0.1in}For any graph,
$\operatorname{ir}_{b}\leq\gamma_{b}\leq\gamma\leq\Gamma\leq\Gamma_{b}%
\leq\operatorname{IR}_{b}.$
\end{corollary}

\section{Properties of Irredundant Broadcasts}

\label{Sec_Ir}

\begin{figure}[ptb]%
\centering
\includegraphics[
height=1.7668in,
width=2.8746in
]%
{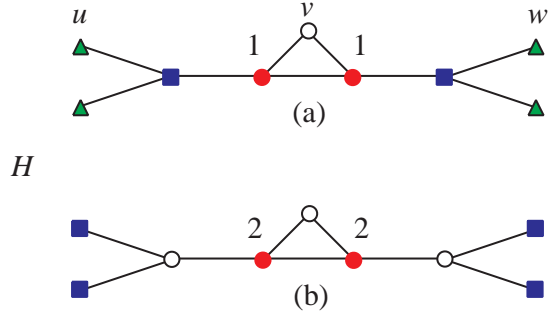}%
\caption{A maximal irredundant set that is not a maximal irredundant
broadcast}%
\label{Fig_ir-set}%
\end{figure}

While a minimal dominating set is also a minimal dominating broadcast, a
maximal irredundant set is not necessarily a maximal irredundant broadcast. In
\cite{Ahmadi} the path $P_{6}$ is used to illustrate this fact. Here we use
the graph $H$ in Figure \ref{Fig_ir-set}. The red vertices (solid circles)
form a maximal irredundant set. Their private neighbours are shown as blue
squares, while the undominated vertices are shown as green triangles. If any
blue, green or white vertex is changed to red, the resulting set of red
vertices is not irredundant. However, the broadcast in Figure \ref{Fig_ir-set}%
(b) is a minimal dominating and a maximal irredundant broadcast (although not
of minimum cost); the vertices in the private boundaries are shown in blue. By
broadcasting with a strength of 3 from $v$ and observing that $\{u,v,w\}$ is a
multipacking, we see that $\operatorname{ir}_{b}(H)=\gamma_{b}%
(H)=3>2=\operatorname{ir}(H)$. The difference $\operatorname{ir}%
_{b}-\operatorname{ir}$ can be arbitrarily large for connected graphs: form
the graph $J_{k}$ from $k$ copies $H_{1},...,H_{k}$ of $H$ by joining a
rightmost (with respect to the representation in Figure \ref{Fig_ir-set})
undominated vertex of $H_{i-1}$ to a leftmost undominated vertex of
$H_{i},\ i=2,...,k$. Then $\operatorname{ir}(J_{k})=2k$ and $\operatorname{ir}%
_{b}(J_{k})=\gamma_{b}(J_{k})=3k$. The ratio $\operatorname{ir}_{b}%
/\operatorname{ir}$, however, is bounded: as shown in \cite{BMT, LauraT},
$\gamma_{b}\leq\frac{3}{2}\operatorname{ir}$ for all graphs, hence
$\operatorname{ir}_{b}/\operatorname{ir}\leq\gamma_{b}/\operatorname{ir}%
\leq\frac{3}{2}$. The graphs $J_{k},\ k\geq1$, show that the bound is tight.

On the other hand, the tree obtained by subdividing each edge of $K_{1,n}$
once satisfies $\operatorname{ir}=\gamma=n$ and $\operatorname{ir}_{b}%
=\gamma_{b}=2$, hence $\operatorname{ir}-\operatorname{ir}_{b}$ and
$\operatorname{ir}/\operatorname{ir}_{b}$ can be arbitrarily large.%

\subsection{Maximal irredundant broadcasts}

In order to find graphs with $\operatorname{ir}_{b}<\gamma_{b}$ and to bound
the ratio $\gamma_{b}/\operatorname{ir}_{b}$ we first need to answer the
following question:

\begin{question}
\label{Ques_max-irr}What makes an irredundant broadcast maximal irredundant?
\end{question}

Let $f$ be a broadcast on a graph $G=(V,E)$. Recall that $U_{f}$ denotes the
set of all vertices not dominated by $f$. Suppose $U_{f}\neq\varnothing$. We
say a vertex $v\in V(G)$ is \emph{blocked} if some shortest path from $v$ to
$U_{f}$ contains a vertex in $V_{f}^{+}-\{v\}$. Define

\begin{itemize}
\item $\beta_{f}=\{v\in V(G):v$ is blocked$\}$

\item $V_{f}^{\ast}=\{v\in V_{f}^{+}:v$ is not blocked$\}=V_{f}^{+}-\beta_{f}$,
\end{itemize}

\noindent and for each $v\in V-U_{f}$,

\begin{itemize}
\item $d_{f}(v)=\min\{d(b,u):b\in B_{f}(v),\ u\in U_{f}\}$, that is,
$d_{f}(v)$ is the minimum distance between the boundary of $v$ and an
undominated vertex, which is also the smallest integer $k$ such that
$U_{f}\cap N_{f(v)+k}[v]\neq\varnothing$,

\item $U_{f}(v)=\{u\in U_{f}:d(u,v)=d_{f}(v)+f(v)\}$.
\end{itemize}

Given a sequence $s=v_{0},...,v_{k}$ of vertices of $G$, define its\emph{
associated sequence} $s_{f}=f_{0},...,f_{k}$ of broadcasts by $f_{0}=f$ and,
for $i=0,...,k-1$,%
\[
f_{i+1}=\left\{
\begin{tabular}
[c]{ll}%
$(f_{i}-\{(v_{i},f_{i}(v_{i}))\})\cup\{(v_{i},f_{i}(v_{i})+d_{f_{i}}%
(v_{i}))\}$ & if $U_{f_{i}}\neq\varnothing$\\
$f_{i}$ & otherwise.
\end{tabular}
\ \right.
\]
If $u\in U_{f}$, $w\in V_{f}^{+}$ and $\operatorname{PB}_{f}(w)\subseteq
N(u)$, we say that $u$ $f$-\emph{annihilates} $w$. Similarly, if $v\in
V-U_{f}$, $w\in V_{f}^{+}$ and $\operatorname{PB}_{f}(w)\subseteq
N_{f(v)+d_{f}(v)}[v]$, we say that $v$ $f$-\emph{annihilates} $w$. If the
broadcast is clear from the context, we simply say that $v$ \emph{annihilates}
$w$.

\begin{theorem}
\label{ThmMax_irr}An irredundant broadcast $f$ is maximal irredundant if and
only if the following conditions hold.

\begin{enumerate}
\item[$(i)$] Each $w\in V-\beta_{f}-V_{f}^{+}$ $f$-annihilates some $v\in
V_{f}^{+}$.

\item[$(ii)$] For each $v_{0}\in V_{f}^{\ast}$ there exist a finite sequence
$s=v_{0},...,v_{t}$ of (not necessarily distinct) vertices in $V_{f}^{+}$ and
its associated sequence $s_{f}=f_{0},...,f_{t}$ of broadcasts such that
$v_{i}$ $f_{i}$-annihilates $v_{i+1}$, $i=0,...,t-1$, and either $v_{t}%
\in\beta_{f_{t}}$ or $U_{f_{t}}=\varnothing$.
\end{enumerate}
\end{theorem}

\begin{figure}[ptb]%
\centering
\includegraphics[
height=2.2105in,
width=5.9499in
]%
{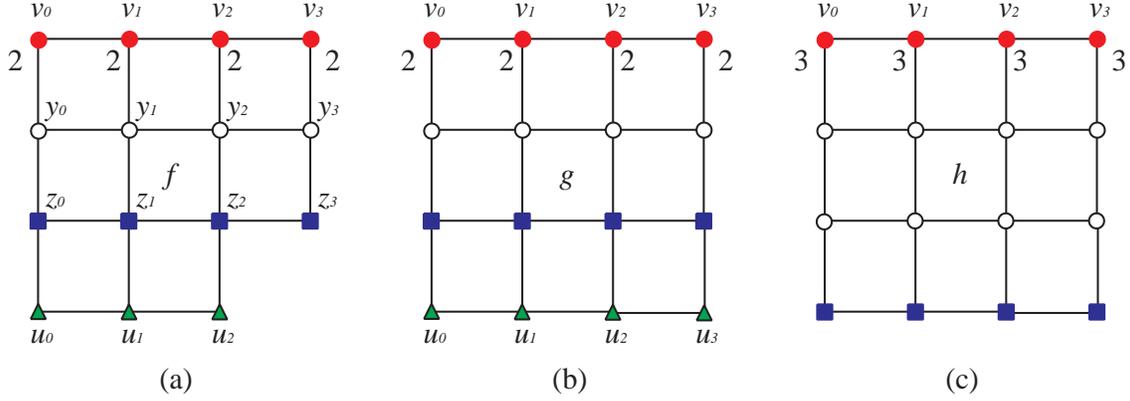}%
\caption{(a) A maximal irredundant broadcast $f$,\ (b) a non-maximal
irredundant broadcast $g$, and (c) a minimal dominating and maximal
irredundant broadcast $h$}%
\label{FigMax_ir}%
\end{figure}

Before proving Theorem \ref{ThmMax_irr} we illustrate $(ii)$ with an example.
Consider the broadcast $f=f_{0}$ in Figure \ref{FigMax_ir}(a), where we denote
the vertices in $V_{f}^{+}$ by red (solid) circles, the vertices that belong
to some private boundary by blue squares, and the vertices in $U_{f}$ by green
triangles. Since each vertex in $V_{f}^{+}=\{v_{0},...,v_{3}\}$ has a nonempty
$f$-private boundary, $f$ is irredundant. Note that $d_{f}(v_{0})=d_{f}%
(v_{1})=d_{f}(v_{2})=1$ and $d_{f}(v_{3})=2$. For $v_{0}\in V_{f}^{\ast}$ and
$f_{1}=(f_{0}-\{(v_{0},2)\})\cup\{(v_{0},3)\}$, $\operatorname{PB}_{f_{1}%
}(v_{1})=\varnothing$, hence $f_{1}$ is not irredundant. Now $d_{f_{1}}%
(v_{1})=1$ and we let $f_{2}=(f_{1}-\{(v_{1},2)\})\cup\{(v_{1},3)\}$. Although
$u_{0}\in\operatorname{PB}_{f_{2}}(v_{0})$ and $u_{1}\in\operatorname{PB}%
_{f_{2}}(v_{1})$, $\operatorname{PB}_{f_{2}}(v_{2})=\varnothing$, hence
$f_{2}$ is not irredundant. Since $d_{f_{2}}(v_{2})=1$, we define
$f_{3}=(f_{2}-\{(v_{2},2)\})\cup\{(v_{2},3)\}$. We see that for $i=0,1,2$,
$u_{i}\in\operatorname{PB}_{f_{3}}(v_{i})$, but $\operatorname{PB}_{f_{3}%
}(v_{3})=\varnothing=U_{f_{3}}$. Therefore, for $v_{0}\in V_{f}^{\ast}$, the
sequences $s=v_{0},...,v_{3}$ and $s_{f}$ satisfy $(ii)$, and $f$ cannot be
extended by \emph{this specific sequence} of broadcasts to a larger maximal
broadcast. The obvious subsequences for $v_{1}$ and $v_{2}$ also satisfy
$(ii)$. Since $v_{3}\in\beta_{f}$, we do not need to find a sequence beginning
with $v_{3}$. Also checking $(i)$ (for $y_{0},...,y_{3},z_{0},...,z_{3}%
,u_{0},...,u_{2}$) we deduce that $f$ is maximal irredundant. Similarly
checking Theorem \ref{ThmMax_irr}$(ii)$ for the irredundant broadcast $g$ in
Figure \ref{FigMax_ir}(b) shows that $g$ is not maximal irredundant, and
indeed the irredundant broadcast $h$ in Figure \ref{FigMax_ir}(c) satisfies
$h>g$.\bigskip%

\noindent\textbf{Proof of Theorem \ref{ThmMax_irr}}.\hspace{0.1in}Suppose $f$
is a maximal irredundant broadcast. If $f$ is dominating, then $U_{f}%
=\varnothing$ and there is nothing to prove, so assume $U_{f}\neq\varnothing$.

\noindent$(i)\hspace{0.1in}$First consider $u\in U_{f}$ and let
$g=(f-\{(u,0)\})\cup\{(u,1)\}$. Then $u\in\operatorname{PB}_{g}(u)$ and thus
$\operatorname{PB}_{g}(u)\neq\varnothing$. Since $g$ is not irredundant,
$\operatorname{PB}_{g}(v)=\varnothing$ for some $v\in V_{f}^{+}$. Therefore
$\operatorname{PB}_{f}(v)\subseteq N_{g}[u]=N[u]$. Since $u\in U_{f}$,
$u\notin\operatorname{PB}_{f}(v)$, thus $\operatorname{PB}_{f}(v)\subseteq
N(u)$ and $(i)$ holds.

Now consider $w\in V-U_{f}$ such that $f(w)=0$ and define $f^{\prime
}=(f-\{(w,0)\})\cup\{(w,d_{f}(w))\}$. Since $U_{f}\neq\varnothing$,
$U_{f}(w)\neq\varnothing$, hence $\operatorname{PB}_{f_{1}}(w)=U_{f}%
(w)\neq\varnothing$. Since $d_{f}(w)>0$, $f^{\prime}>f$. By the maximality of
$f$, $\operatorname{PB}_{f_{1}}(v)=\varnothing$ for some $v\in V_{f}^{+}$,
which implies that $\operatorname{PB}_{f}(v)\subseteq N_{d_{f}(w)}[w]$. Hence
$w$ annihilates $v$, as required.

\noindent$(ii)\hspace{0.1in}$Let $f_{0}=f$ and consider any $v_{0}\in
V_{f}^{\ast}$ and the broadcast $f_{1}=(f_{0}-\{(v_{0},f_{0}(v_{0}%
))\})\cup\{(v_{0},f_{0}(v_{0})+d_{f}(v_{0}))\}$. Let $u_{0}\in U_{f}(v_{0})$.
Since $d(u_{0},v_{0})=f_{1}(v_{0})$ and $d(u_{0},v^{\prime})>f(v^{\prime
})=f_{1}(v^{\prime})$ for all $v^{\prime}\in V_{f}^{+}-\{v_{0}\}$, $u_{0}%
\in\operatorname{PB}_{f_{1}}(v_{0})$. If $\operatorname{PB}_{f_{1}}(v^{\prime
})\neq\varnothing$ for all $v^{\prime}\in V_{f}^{+}-\{v_{0}\}$, then $f_{1}$
is irredundant, contradicting the maximality of $f$. Hence $\operatorname{PB}%
_{f_{1}}(v_{1})=\varnothing$ for some $v_{1}\in V_{f}^{+}-\{v_{0}\}$. If also
$v_{1}\in V_{f}^{+}-V_{f_{1}}^{\ast}$ or $U_{f_{1}}=\varnothing$, we are done,
hence assume $v_{1}\in V_{f_{1}}^{\ast}$ and $U_{f_{1}}\neq\varnothing$.
Define $f_{2}=(f_{1}-\{(v_{1},f_{1}(v_{1}))\})\cup\{(v_{1},f_{1}%
(v_{1})+d_{f_{1}}(v_{1}))\}$ and let $u_{1}\in U_{f_{1}}(v_{1})$. Then
$u_{1}\in\operatorname{PB}_{f_{2}}(v_{1})$. Continue as before. Since each
vertex has finite eccentricity and $f_{i}(v)\leq e(v)$ for each $i$ and each
vertex $v$, the process ends. Hence eventually we obtain sequences
$s=v_{0},...,v_{t}$ and $s_{f}=f_{0},...,f_{t}$ such that $\operatorname{PB}%
_{f_{i}}(v_{i})=\varnothing$ for $i=1,...,t$ and $v_{t}\in V_{f}^{+}-V_{f_{t}%
}^{\ast}$ or $U_{f_{t}}=\varnothing$.

Conversely, suppose $f$ is an irredundant broadcast such that $(i)$ and $(ii)$
hold, but $f$ is not maximal irredundant. Let $g>f$ be an irredundant
broadcast. We first prove a lemma.

\bigskip

\noindent\textbf{Lemma \ref{ThmMax_irr}.1\hspace{0.1in}}There exists $v\in
V_{f}^{\ast}$ such that $g(v)>f(v)$.

\noindent\textbf{Proof of Lemma \ref{ThmMax_irr}.1\hspace{0.1in}}Let $x$ be a
vertex such that $g(x)>f(x)$. If $x\in U_{f}$, then, by $(i)$, there exists
$y\in V_{f}^{+}$ such that $\operatorname{PB}_{f}(y)\subseteq N(x)$. Since $g$
is irredundant, $\operatorname{PB}_{g}(y)\neq\varnothing$. Evidently, then,
$g(y)>f(y)$.

Now assume $x\in V-U_{f}$, $f(x)=0$ and $g(x)>f(x)$. Since $g$ is irredundant,
$\operatorname{PB}_{g}(x)\neq\varnothing$. Therefore some $u\in U_{f}$ belongs
to $\operatorname{PB}_{g}(x)$. Now $g(x)=d(x,u)\geq d_{f}(x)$, and by $(i)$,
$x$ $f$-annihilates some $y\in V_{f}^{+}$, that is, $\operatorname{PB}%
_{f}(y)\subseteq N_{d_{f}(x)}[x]$. Therefore $\operatorname{PB}_{f}%
(y)\subseteq N_{g}[x]$. Since $g$ is irredundant, $\operatorname{PB}%
_{g}(y)\neq\varnothing$, which implies that $g(y)>f(y)$.

If $y\in V_{f}^{\ast}$, we are done, hence assume $y$ is blocked. Then some
shortest path from $y$ to $U_{f}$ contains a vertex in$\ V_{f}^{+}-\{y\}$.
Amongst all such paths, let $Q$ be one that contains the most vertices in
$V_{f}^{+}$. Let $u\in U_{f}$ be the terminal vertex of $Q$ and let $v$ be the
vertex in $V_{f}^{+}$ nearest to $u$. Since $Q$ is a shortest $y-u$ path
containing $v$, the $v-u$ subpath $P$ of $Q$ is a shortest path from $v$ to
$U_{f}$. By the choice of $v$, $v$ is the only vertex in $V_{f}^{+}$ on $P$.
If any shortest path $P^{\prime}$ from $v$ to $U_{f}$ contains a vertex from
$V_{f}^{+}-\{v\}$, then $(Q-P)\cup P^{\prime}$ is a shortest path from $y$ to
$U_{f}$ that contains more vertices from $V_{f}^{+}$ than $Q$ does, contrary
to the choice of $Q$. Hence $v\in V_{f}^{\ast}$. Since $g$ is irredundant,
$\operatorname{PB}_{g}(y)\neq\varnothing$, hence some vertex $u^{\prime}\in
U_{f}$ belongs to $\operatorname{PB}_{g}(y)$. Therefore $g(y)\geq
f(y)+d_{f}(y)$ and $u$ hears the broadcast $g$ from $y$ (and $u\in
\operatorname{PB}_{g}(y)$ if and only if $g(y)=f(y)+d_{f}(y)$). Since $v$ lies
on $Q$, $N_{f}[v]\cup U_{f}[v]\subseteq N_{g}[y]$ and, in particular,
$\operatorname{PB}_{f}(v)\subseteq N_{g}[y]$. But then, since $g$ is
irredundant, $g(v)>f(v)$ to ensure that $\operatorname{PB}_{g}(v)\neq
\varnothing$.$~\blacklozenge$

\bigskip

Consider $v\in V_{f}^{\ast}$ such that $g(v)>f(v)$. Let $v_{0}=v$ and
$u_{0}\in\operatorname{PB}_{g}(v_{0})$; note that $u_{0}\in U_{f}$ and
$f(v_{0})+d_{f}(v_{0})\leq g(v_{0})$. Let $s=v_{0},...,v_{t}$ and $s_{f}%
=f_{0},...,f_{t}$ be the sequences guaranteed by $(ii)$. By definition,
$f_{1}=(f-\{(v_{0},f(v_{0}))\})\cup\{(v_{0},f(v_{0})+d_{f}(v_{0}))\}$ and
$\operatorname{PB}_{f_{1}}(v_{1})=\varnothing$, and therefore
$\operatorname{PB}_{f_{0}}(v_{1})\subseteq N_{f_{1}}[v_{0}]$. Since
$g(v_{0})\geq f(v_{0})+d_{f}(v_{0})=f_{1}(v_{0})$, $\operatorname{PB}_{f_{0}%
}(v_{1})\subseteq N_{g}[v_{0}]$. Since $g$ is irredundant, $\operatorname{PB}%
_{g}(v_{1})\neq\varnothing$. Hence $g(v_{1})>f_{0}(v_{1})=f_{1}(v_{1})$ and
there exists $u_{1}\in\operatorname{PB}_{g}(v_{1})\cap U_{f_{1}}$. As in the
case of $u_{0},\ v_{0}$ and $f$, we see that $g(v_{1})\geq f_{1}%
(v_{1})+d_{f_{1}}(v_{1})=f_{2}(v_{1})$, so that $\operatorname{PB}_{f_{1}%
}(v_{2})\subseteq N_{f_{2}}[v_{1}]\subseteq N_{g}[v_{1}]$. Continuing the
process we see that $g(v_{i})\geq f_{i+1}(v_{i})$, $i=1,...,t-1$ and thus
$\operatorname{PB}_{f_{t-1}}(v_{t})\subseteq N_{f_{t}}[v_{t-1}]\subseteq
N_{g}[v_{t-1}]$. Since $g$ is irredundant, $\operatorname{PB}_{g}(v_{t}%
)\neq\varnothing$, hence there exists $u_{t}\in\operatorname{PB}_{g}%
(v_{t})\cap U_{f_{t}}$. Then $U_{f_{t}}\neq\varnothing$ and so, by $(ii)$,
$v_{t}\in V_{f}^{+}-V_{f_{t}}^{\ast}$.

Now, since $\operatorname{PB}_{f_{t-1}}(v_{t})\subseteq N_{g}[v_{t-1}]$ and
$\operatorname{PB}_{g}(v_{t})\neq\varnothing$, $g(v_{t})>f_{t-1}(v_{t})\geq
f(v_{t})$. As in the proof of Lemma \ref{ThmMax_irr}.1, there exists
$v^{\prime}\in V_{f_{t}}^{\ast}$ such that $g(v^{\prime})>f_{t}(v^{\prime})$.
Let $v_{0}^{\prime}=v^{\prime}$, $f^{\prime}=f_{t}$ and $u_{0}^{\prime}%
\in\operatorname{PB}_{g}(v_{0}^{\prime})$.Let $s^{\prime}=v_{0}^{\prime
},...,v_{t^{\prime}}^{\prime}$ and $s_{f^{\prime}}^{\prime}=f_{0}^{\prime
},...,f_{t^{\prime}}^{\prime}$ be the sequences for $v^{\prime}$ and
$f^{\prime}$ guaranteed by $(ii)$. We repeat the above argument until we
obtain $u_{t^{\prime}}^{\prime}\in\operatorname{PB}_{g}(v_{t^{\prime}}%
^{\prime})\cap U_{f_{t^{\prime}}^{\prime}}$, so $U_{f_{t^{\prime}}^{\prime}%
}\neq\varnothing$ and thus, by $(ii)$, $v_{t^{\prime}}^{\prime}\in
V_{f^{\prime}}^{+}-V_{f_{t^{\prime}}^{\prime}}^{\ast}$. However, since
$|U_{f}|$ is finite, this process cannot continue indefinitely, and eventually
we obtain sequences $s^{(k)}=v_{0}^{(k)},...,v_{t^{(k)}}^{(k)}$ and
$s_{f^{(k)}}^{(k)}=f_{0}^{(k)},...,f_{t^{(k)}}^{(k)}=h$ such that
$U_{h}=\varnothing$, the final contradiction. Hence $f$ is maximal
irredundant, as required.~$\blacksquare$

\begin{corollary}
If $f$ is a maximal irredundant broadcast, then each vertex in $U_{f}$ is at
distance $f(v)+1$ from some vertex $v\in V_{f}^{+}$.
\end{corollary}

\noindent\textbf{Proof.\hspace{0.1in}}If $u\in U_{f}$ and $d(u,v)\geq f(v)+2$
for all $v\in V_{f}^{+}$, let $g=(f-\{(u,0)\})\cup\{(u,1)\}$. Then
$u\in\operatorname{PB}_{g}(u)$ and thus $\operatorname{PB}_{g}(u)\neq
\varnothing$. Also, $\operatorname{PB}_{g}(v)=\operatorname{PB}_{f}%
(v)\neq\varnothing$ for all $v\in V_{f}^{+}$, which is a
contradiction.~$\blacksquare$

\begin{corollary}
If an irredundant broadcast is dominating, then it is maximal irredundant and
minimal dominating.
\end{corollary}

\noindent\textbf{Proof.\hspace{0.1in}}Let $f$ be a dominating irredundant
broadcast. Then $U_{f}=\varnothing$ and $(i)$ and $(ii)$ of Theorem
\ref{ThmMax_irr} are vacuously true. Hence $f$ is maximal irredundant. Also,
$f$ is minimal dominating by the definition of broadcast
irredundance.~$\blacksquare$

\section{Comparing broadcast domination and irredundance numbers}

\label{Sec_Bound}As mentioned in Section \ref{Sec_Ir}, $\operatorname{ir}%
_{b}\leq\gamma_{b}(G)\leq\frac{3}{2}\operatorname{ir}(G)$ for all graphs $G$,
the second inequality being established in \cite{BMT}. Although
$\operatorname{ir}$ and $\operatorname{ir}_{b}$ are not comparable, it is
reasonable to expect that $\gamma_{b}(G)\leq(1+c)\operatorname{ir}_{b}(G)$ for
some constant $c$. That this is indeed the case is the subject of our main result.

\bigskip\label{here}

\noindent\textbf{Theorem \ref{Thm_irb-gammab_bound}\hspace{0.1in}}\emph{For
any graph }$G$\emph{, }$\gamma_{b}(G)\leq\frac{5}{4}\operatorname{ir}_{b}%
(G)$\emph{.}

\bigskip

\noindent\textbf{Proof.\hspace{0.1in}}Let $f$ be any maximal irredundant
broadcast on $G$. We construct a dominating broadcast $g$ of $G$ with
$\sigma(g)\leq\frac{5}{4}\sigma(f)$. Define the graph $H$ as follows:
$V(H)=V_{f}^{+}$, and vertices $u$ and $v$ of $H$ are adjacent if and only if
$N_{f}[u]\cap N_{f}[v]\neq\varnothing$. Let $X_{i},\ i=1,...,t$, denote the
(vertex sets of the) components of $H$. Define the collections $\mathcal{A}$
and $\mathcal{A}_{1}-\mathcal{A}_{4}$ of components of $H$ by%
\begin{align*}
\mathcal{A}_{i}  &  =\{X_{j}:%
{\textstyle\sum_{x\in X_{j}}}
f(x)=j,\ j=1,2,3\}\\
\mathcal{A}_{4}  &  =\{X_{j}:%
{\textstyle\sum_{x\in X_{j}}}
f(x)\geq4\}\\
\mathcal{A}  &  =%
{\textstyle\bigcup_{i=1}^{4}}
\mathcal{A}_{i}.
\end{align*}
For $X\in\mathcal{A}$, let

\begin{itemize}
\item[\ ] $G_{X}^{\prime}$ be the subgraph of $G$ induced by $%
{\textstyle\bigcup_{x\in X}}
N_{f}[x]$ and let

\item[\ ] $G_{X}$ be the subgraph of $G$ induced by $%
{\textstyle\bigcup_{x\in X}}
N_{f}[x]$ together with all vertices in $U_{f}$ that annihilate a vertex in
$X$.
\end{itemize}

\noindent We first establish a bound on the diameters and radii of $G_{X}$ and
$G_{X}^{\prime}$.

\smallskip

\noindent\textbf{Lemma \ref{Thm_irb-gammab_bound}.1\hspace{0.1in}}For each
$X\in\mathcal{A},\ \operatorname{diam}(G_{X})\leq2%
{\textstyle\sum_{x\in X}}
f(x)+2,\ \operatorname{rad}(G_{X})\leq%
{\textstyle\sum_{x\in X}}
f(x)+1,\ \operatorname{diam}(G_{X}^{\prime})\leq2%
{\textstyle\sum_{x\in X}}
f(x)$ and $\operatorname{rad}(G_{X}^{\prime})\leq%
{\textstyle\sum_{x\in X}}
f(x).\smallskip$

\noindent\textbf{Proof of Lemma \ref{Thm_irb-gammab_bound}.1.\hspace{0.1in}%
}Say $X=\{a_{1},...,a_{t}\}$ and let $P_{X}$ be a diametrical path of $G_{X}%
$.\textbf{ }The diameter of $G_{X}$ is maximized when $(i)$ all $a_{i}$ lie on
$P_{X}$, $(ii)$ if (without loss of generality) $a_{1}$ is the first vertex of
$X$ on $P_{X}$, then $a_{1}$ is preceded on $P_{X}$ by $f(a_{1})$ vertices in
$N_{f}[a_{1}]$ and a vertex $u_{1}\in U_{f}$, which is the origin of $P_{X}$,
$(iii)$ a similar comment holds for the last vertex, say $a_{t}$, of $X$ on
$P_{X}$, where some vertex $u_{t}\in U_{f}$ is the terminus of $P_{X}$, and
$(iv)$ for any $a_{i}$, if $a_{i+1}$ is the next vertex of $X$ on $P_{X}$,
then $d(a_{i},a_{i+1})=f(a_{i})+f(a_{i+1})$. Hence
\begin{align*}
\operatorname{diam}(G_{X})  &  \leq\lbrack1+f(a_{1})]+[f(a_{1})+f(a_{2}%
)]+\cdots+[f(a_{t-1})+f(a_{t})]+[f(a_{t})+1]\\
&  =2%
{\textstyle\sum_{i=1}^{t}}
f(a_{i})+2.
\end{align*}
Moreover, if at least one of $u_{1}$ and $u_{t}$ exists, then any central
vertex of $P_{X}$ has eccentricity at most $%
{\textstyle\sum_{x\in X}}
f(x)+1$, hence $\operatorname{rad}(G_{X})\leq%
{\textstyle\sum_{x\in X}}
f(x)+1$. The argument for $\operatorname{diam}(G_{X}^{\prime})$ and
$\operatorname{rad}(G_{X}^{\prime})$ is the same, except that vertices in
$N_{f}[a_{1}]$ and $N_{f}[a_{t}]$ are the origin and terminus of a diametrical
path.$~\blacklozenge$

\bigskip

\begin{figure}[ptb]%
\centering
\includegraphics[
height=2.1119in,
width=4.2557in
]%
{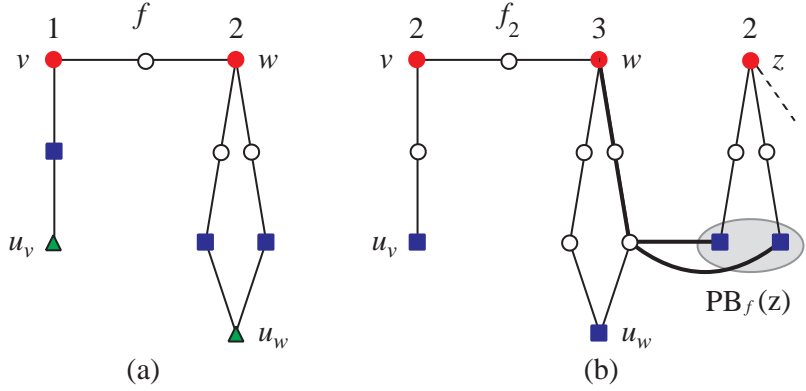}%
\caption{Irredundant and minimal dominating broadcasts on $G_{X}$, where
$X=\{v,w\}\in\mathcal{B}$}%
\label{Fig_B}%
\end{figure}

Suppose $X=\{v,w\}$, where $f(v)=1$, $f(w)\in\{1,2\}$, $d(v,w)=f(w)$,
$u_{v}\in U_{f}$ annihilates $v$ and is nonadjacent to all vertices in
$B_{f}(w)$, while $u_{w}\in U_{f}$ annihilates $w$ and is nonadjacent to all
vertices in $B_{f}(v)$. See Figure \ref{Fig_B}(a). Then $d_{f}(v)=d_{f}(w)=1$.
Since $v$ hears the broadcast from $w$, $\operatorname{diam}(G_{X}%
)\leq2[f(v)+f(w)]+1$. Consider the broadcasts%
\begin{align*}
f_{1}  &  =(f-\{(v,1)\})\cup\{(v,2)\}\text{ and}\\
f_{2}  &  =(f_{1}-\{(w,f(w))\})\cup\{(w,f(w)+1)\}.
\end{align*}
Note that $u_{v}\in\operatorname{PB}_{f_{1}}(v)\cap\operatorname{PB}_{f_{2}%
}(v)$ and $u_{w}\in\operatorname{PB}_{f_{2}}(w)$ (see Figure \ref{Fig_B}(b)).
By the maximality of $f$ there exists a vertex $z\in V_{f}^{+}-X$ such that
$\operatorname{PB}_{f}(z)\subseteq B_{f_{2}}(v)\cup B_{f_{2}}(w)$; we also say
that $X$ \emph{annihilates} $z$. Let

\begin{itemize}
\item[\ ] $\mathcal{B}$ be the subset of $\mathcal{A}_{2}\cup\mathcal{A}_{3}$
that consists of all such $X$.
\end{itemize}

\smallskip%

\noindent\textbf{Lemma \ref{Thm_irb-gammab_bound}.2\hspace{0.1in}}Each
$X\in\mathcal{B}$ annihilates some $z\in X^{\prime}\in\mathcal{A}%
-\mathcal{A}_{1}$.

$\smallskip$

\noindent\textbf{Proof of Lemma \ref{Thm_irb-gammab_bound}.2.\hspace{0.1in}%
}Suppose $X\in\mathcal{B}$ annihilates $z\in X^{\prime},\ X^{\prime}%
\in\mathcal{A}_{1}$. Then $X^{\prime}=\{z\}$, $f(z)=1$ and, by definition of
$H$, no vertex in $N[z]$ hears $f$ from a vertex in $V_{f}^{+}-\{z\}$. Hence
$\operatorname{PB}_{f}(z)=N[z]$. Since $X$ annihilates $z$, each vertex in
$\operatorname{PB}_{f}(z)$ is adjacent to a vertex in $B_{f}(x)$, for some
$x\in X$. But then $z$ itself is adjacent to a vertex $b\in B_{f}(x)$, for
some $x\in X$, so $b\in B_{f}(z)\cap B_{f}(x)$, contrary to $z$ being an
isolated vertex of $H$.$~\blacklozenge$

\bigskip

Suppose $X=\{v,w\}$, where $f(v)=1$, $f(w)=2$, $d(v,w)=3$, $u_{w}\in U_{f}$
annihilates $w$, and some vertex $v^{\prime}\in\operatorname{PB}_{f}(v)$ is
nonadjacent to all vertices in $B_{f}(w)$. Then $d_{f}(w)=1$, and since
$v^{\prime}$ is nonadjacent to $B_{f}(w)$, $w$ does not annihilate $v$. Since
$v$ does not hear $f$ from any other vertex, $v\in\operatorname{PB}_{f}(v)$.
Hence no vertex in $U_{f}$ annihilates $v$ either. Therefore
$\operatorname{diam}(G_{X})\leq2[f(v)+f(w)]+1=7$. Since $w$ does not
annihilate $v$, the maximality of $f$ implies that $w$ annihilates a vertex
$z\in X^{\prime}\neq X$. Let

\begin{itemize}
\item[\ ] $\mathcal{C}$ be the subset of $\mathcal{A}_{3}$ that consists of
all such $X$.
\end{itemize}

The proof of the next lemma is similar to that of Lemma
\ref{Thm_irb-gammab_bound}.2 and is omitted.

\smallskip

\noindent\textbf{Lemma \ref{Thm_irb-gammab_bound}.3\hspace{0.1in}}For each
$X\in\mathcal{C}$, the vertex $w\in X$ such that $f(w)=2$ annihilates some
$z\in X^{\prime}\in\mathcal{A}-\mathcal{A}_{1}$.

\medskip

\begin{figure}[ptb]%
\centering
\includegraphics[
height=2.5028in,
width=4.3673in
]%
{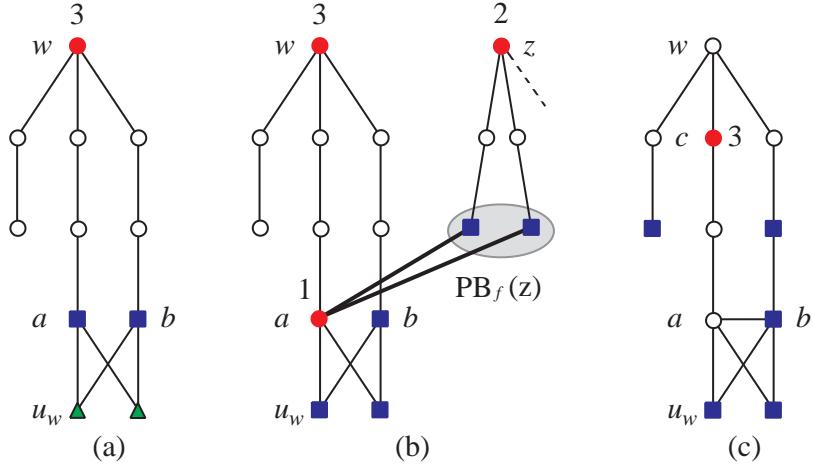}%
\caption{Irredundant and minimal dominating broadcasts on $G_{X}$, where
$X=\{w\}\in\mathcal{D}$}%
\label{Fig_D}%
\end{figure}

Suppose $X=\{w\}$, where $f(w)\in\{2,3\}$, $u_{w}\in U_{f}$ annihilates $w$
and $\operatorname{PB}_{f}(w)$ contains two nonadjacent vertices $a$ and $b$.
Let $f_{1}=(f-\{(a,0)\}\cup\{(a,1)\}$. Then $u_{w}\in\operatorname{PB}_{f_{1}%
}(a)$ and $b\in\operatorname{PB}_{f_{1}}(w)$. See Figure \ref{Fig_D}(a) and
(b). By the maximality of $f$ there exists a vertex $z\in V_{f}^{+}-\{w\}$
such that $a$ annihilates $z$, that is, $\operatorname{PB}_{f}(z)\subseteq
N[a]$. Let $f_{2}=(f-\{(w,f(w))\}\cup\{(w,f(w)+1)\}$. Then $\operatorname{PB}%
_{f}(z)\subseteq N_{f_{2}}[w]$, hence $w$ also annihilates $z$. Let

\begin{itemize}
\item[\ ] $\mathcal{D}$ be the subset of $\mathcal{A}_{2}\cup\mathcal{A}_{3}$
that consists of all such $X$, and

\item[\ ] $\mathcal{E}=\mathcal{B}\cup\mathcal{C}\cup\mathcal{D}$.
\end{itemize}

The proof of the next lemma is also omitted.

\smallskip

\noindent\textbf{Lemma \ref{Thm_irb-gammab_bound}.4\hspace{0.1in}}Each
$X\in\mathcal{D}$ annihilates some $z\in X^{\prime}\in\mathcal{A}%
-\mathcal{A}_{1}$.

\medskip%

Consider $X\in\mathcal{E}$; without loss of generality, say $X=X_{1}$, and
that $X_{1}$ annihilates $z\in X_{2}\in\mathcal{A}-\mathcal{A}_{1}$. (It is
possible that $X_{1}$ annihilates vertices of several sets $X^{\prime}$, in
which case we choose an arbitrary $X^{\prime}=X_{2}$.) Let%
\[
G_{1,2}\text{ be the subgraph of }G\text{ induced by }V(G_{X_{1}})\cup
V(G_{X_{2}}).
\]
If $X_{1}\in\mathcal{B}$, both the origin and terminal of the diametrical path
$P_{X_{1}}$ belong to $U_{f}$. Since $X_{1}$ annihilates $z$, one of these
ends of $P_{X_{1}}$ does not belong to any diametrical path of $G_{1,2}$. If
$X_{1}\in\mathcal{C}\cup\mathcal{D}$, then $w$ with $f(w)\in\{2,3\}$ and
$u_{w}\in U_{f}(w)$ annihilates $z$, hence no vertex in $U_{f}(w)$ belongs to
any diametrical path of $G_{1,2}$. Similarly, if some vertex $u_{z}\in
U_{f}(z)$ belongs to a diametrical path $P_{X_{2}}$ of $X_{2}$, then $u_{z}$
does not belong to a diametrical path of $G_{1,2}$. Hence
\begin{align*}
\operatorname{diam}(G_{1,2})  &  \leq\lbrack2%
{\textstyle\sum_{x\in X_{1}}}
f(x)+1]+1+[2%
{\textstyle\sum_{x\in X_{2}}}
f(x)+2]-2\\
&  =2%
{\textstyle\sum_{x\in X_{1}\cup X_{2}}}
f(x)+2.
\end{align*}
Similar to the proof of Lemma \ref{Thm_irb-gammab_bound}.1,%
\[
\operatorname{rad}(G_{1,2})\leq%
{\textstyle\sum_{x\in X_{1}\cup X_{2}}}
f(x)+1.
\]
If there exists $X_{3}\in\mathcal{E}-\{X_{1},X_{2}\}$ that annihilates a
vertex $z^{\prime}$ of $X_{1}$ or $X_{2}$, let%
\[
G_{1,2,3}\text{ be the subgraph of }G\text{ induced by }V(G_{X_{1}})\cup
V(G_{X_{2}})\cup V(G_{X_{3}}).
\]
The diameter of $G_{1,2,3}$ is maximized if (a) all vertices in $X_{1}\cup
X_{2}\cup X_{3}$ lie on a diametrical path $P_{1,2,3}$ of $G_{1,2,3}$, (b)
(without loss of generality) the last vertex of $G_{X_{1}}^{\prime}$ on
$P_{X_{1}}$ is adjacent to the first vertex of $G_{X_{2}}^{\prime}$ on
$P_{X_{2}}$, and (c) the first vertex of $G_{X_{3}}^{\prime}$ on $P_{X_{3}}$
is adjacent to the last vertex of $G_{X_{2}}^{\prime}$ on $P_{X_{2}}$ or the
first vertex of $G_{X_{1}}^{\prime}$ on $P_{X_{1}}$. Hence%
\begin{align*}
\operatorname{diam}(G_{1,2,3})  &  \leq\lbrack2%
{\textstyle\sum_{x\in X_{1}}}
f(x)+1]+1+[2%
{\textstyle\sum_{x\in X_{2}}}
f(x)+2]+1+[2%
{\textstyle\sum_{x\in X_{3}}}
f(x)+1]-4\\
&  =2%
{\textstyle\sum_{x\in X_{1}\cup X_{2}\cup X_{3}}}
f(x)+2\text{ and}\\
\operatorname{rad}(G_{1,2,3})  &  \leq%
{\textstyle\sum_{x\in X_{1}\cup X_{2}\cup X_{3}}}
f(x)+1.
\end{align*}
Continue this process to construct the graph $G(k_{1})=G_{1,...,k_{1}}$ with
diameter $2%
{\textstyle\sum_{x\in X_{1}\cup\cdots\cup X_{k_{1}}}}
f(x)+2$ until no $X\in\mathcal{E}-(X_{1}\cup\cdots\cup X_{k_{1}})$ annihilates
a vertex of $X_{1}\cup\cdots\cup X_{k_{1}}$. If $\mathcal{E}-(X_{1}\cup
\cdots\cup X_{k_{1}})\neq\varnothing$, then there exists $X_{k_{1}+1}%
\in\mathcal{E}-(X_{1}\cup\cdots\cup X_{k_{1}})$ that annihilates some
$X_{k_{1}+2}\in\mathcal{A-A}_{1}$. Construct the graph $G(k_{2})=G_{k_{1}%
+1,...,k_{2}}$ with diameter $2%
{\textstyle\sum_{x\in X_{k_{1}+1}\cup\cdots\cup X_{k_{2}}}}
f(x)+2$ until no $X\in\mathcal{E}-(X_{1}\cup\cdots\cup X_{k_{2}})$ annihilates
a vertex of $X_{1}\cup\cdots\cup X_{k_{2}}$. In this way construct
$G(k_{1}),...,G(k_{s})$ until each $X\in\mathcal{E}$ belongs to some
$G(k_{i})$; that this is possible follows from Lemmas
\ref{Thm_irb-gammab_bound}.2 -- \ref{Thm_irb-gammab_bound}.4. To emphasize,%
\begin{equation}
\operatorname{diam}(G(k_{i}))\leq2%
{\textstyle\sum_{x\in V(G(k_{i}))}}
f(x)+2\text{ and }\operatorname{rad}(G(k_{i}))\leq%
{\textstyle\sum_{x\in V(G(k_{i}))}}
f(x)+1 \label{Eq_Gki}%
\end{equation}
for each $i=1,...,s$. Let
\[
\mathcal{A}_{4}^{\prime}=\{X\in\mathcal{A}_{4}:X\text{ does not belong to any
}G(k_{i})\}.
\]

We now consider the diameters and radii of the graphs $G_{X}$, where
$X\in(\mathcal{A}_{1}\cup\mathcal{A}_{2}\cup\mathcal{A}_{3})-\mathcal{E}$,
more carefully.

\begin{description}
\item[\underline{$X\in\mathcal{A}_{1}$:}] Then $X=\{x\}$, where $f(x)=1$,
$N[x]=\operatorname{PB}_{f}(x)$, $G_{X}=G_{X}^{\prime}$ and
$\operatorname{rad}(G_{X})=1$.

\item[\underline{$X=\{x\}\in\mathcal{A}_{2}\cup\mathcal{A}_{3}$:}]
\label{A2}Then $f(x)\in\{2,3\}$ and $x$ is an isolated vertex of $H$, so no
vertex in $N_{f}[x]$ hears $f$ from any other vertex in $V_{f}^{+}$, and
$B_{f}(x)=\operatorname{PB}_{f}(x)$. Each $u\in U_{f}$ that annihilates $x$ is
adjacent to all vertices in $\operatorname{PB}_{f}(x)$, and, since
$X\notin\mathcal{D}$, $G[\operatorname{PB}_{f}(x)]$ is complete. It follows
that $\operatorname{diam}(G_{X})\leq1+f(x)+[f(x)-1]=2f(x)$. In particular, any
vertex $c$ adjacent to $x$ on an $x-u$ geodesic is a central vertex of $G_{X}$
with eccentricity $e_{G_{X}}(c)=f(x)$. See Figure \ref{Fig_D}(c).

\item[\underline{$X=\{x,y\}\in\mathcal{A}_{2}$:}] Then $f(x)=f(y)=1$. Suppose
$xy\in E(G)$. Since $X\notin\mathcal{B}$, without loss of generality each
vertex in $U_{f}$ that annihilates $y$ is also adjacent to a vertex in
$B_{f}(x)$. Hence $\operatorname{diam}(G_{X})\leq4=2[f(x)+f(y)]$; moreover,
$e_{G_{X}}(x)=2$. If $xy\notin E(G)$, then $d_{G}(x,y)=2=f(x)+f(y)$. Then
neither vertex hears $f$ from the other, so $x\in\operatorname{PB}_{f}(x)$ and
$y\in\operatorname{PB}_{f}(y)$, which means that no vertex in $U_{f}$
annihilates $x$ or $y$. Hence $G_{X}^{\prime}=G_{X}$ and $\operatorname{diam}%
(G_{X})\leq2[f(x)+f(y)]=4$. In addition, any vertex in $N(x)\cap
N(y)\neq\varnothing$ has eccentricity $2$.

\item[\underline{$X=\{x,y\}\in\mathcal{A}_{3}$:}] \label{A3}Without loss of
generality say $f(x)=1$ and $f(y)=2$. Note that $d_{G}(x,y)\geq2$, otherwise
$\operatorname{PB}_{f}(x)=\varnothing$.

Suppose $d_{G}(x,y)=2$. Since $X\notin\mathcal{B}$, at most one of $x$ and $y$
is annihilated by a vertex in $U_{f}$ that is nonadjacent to all vertices in
the boundary of the other. Hence $\operatorname{diam}(G_{X})\leq
f(x)+2f(y)+1=6=2[f(x)+f(y)]$, and $\min\{e(x),e(y)\}=3$.

Suppose $d_{G}(x,y)>2$. Since $x$ and $y$ belong to the same component of $H$,
$d_{G}(x,y)=3$. Moreover, $x\in\operatorname{PB}_{f}(x)$, so no vertex in
$U_{f}$ annihilates $x$. If no vertex in $U_{f}$ annihilates $y$, then
$\operatorname{diam}(G_{X})\leq6=2[f(x)+f(y)]$ and $\operatorname{rad}%
(G_{X})\leq3$.

If some vertex in $U_{f}$ annihilates $y$, then all such vertices in $U_{f}$
are adjacent to all vertices in $\operatorname{PB}_{f}(y)$. Since
$X\notin\mathcal{C}$, each vertex in $\operatorname{PB}_{f}(x)$ is adjacent to
a vertex in $B_{f}(y)$. Therefore $y$ has eccentricity $e_{G_{X}}(y)\leq3$.

\item[\underline{$X=\{x,y,z\}\in\mathcal{A}_{3}$:}] Then $f(x)=f(y)=f(z)=1$.
The diameter of $G_{X}$ is maximized if $x,y$ and $z$ all occur on a
diametrical path and, without loss of generality, $N_{f}[x]\cap N_{f}%
[y]\neq\varnothing$, $N_{f}[y]\cap N_{f}[z]\neq\varnothing$, and $N_{f}[x]\cap
N_{f}[y]=\varnothing$.

Suppose first that $d_{G}(x,y)=d_{G}(y,z)=1$. Then $\operatorname{diam}%
(G_{X})\leq6=2[f(x)+f(y)+f(z)]$.

Now suppose without loss of generality that $d_{G}(x,y)=2$ and $d_{G}(y,z)=1$.
Then $x\in\operatorname{PB}_{f}(x)$, hence no vertex in $U_{f}$ annihilates
$x$, and $\operatorname{diam}(G_{X})\leq6=2[f(x)+f(y)+f(z)]$.

Finally, suppose $d_{G}(x,y)=d_{G}(y,z)=2$. Then $x\in\operatorname{PB}%
_{f}(x)$ and $z\in\operatorname{PB}_{f}(z)$, hence no vertex in $U_{f}$
annihilates $x$ or $z$, and $\operatorname{diam}(G_{X})\leq
6=2[f(x)+f(y)+f(z)]$.

\item In all three cases $e(y)\leq3$, hence $\operatorname{rad}(G_{X})\leq3$.
\end{description}

Our final step is to define a dominating broadcast $g$ with $\sigma
(g)\leq\frac{5}{4}\sigma(f)$. To this end, if $J$ is any subgraph of $G$, we
denote the restriction of $f$ to $J$ by $f\upharpoonleft J$.

\begin{enumerate}
\item Suppose $X\in\mathcal{A}_{4}^{\prime}$, let $c_{X}$ be a central vertex
of $G_{X}$, and define the broadcast $g_{X}$ on $G_{X}$ by%
\[
g_{X}(a)=\left\{
\begin{tabular}
[c]{ll}%
$1+%
{\textstyle\sum_{x\in X}}
f(x)$ & if $a=c_{X}$\\
$0$ & otherwise.
\end{tabular}
\ \ \right.
\]
By Lemma \ref{Thm_irb-gammab_bound}.1, $g_{X}$ dominates $G_{X}$ and, since $%
{\textstyle\sum_{x\in X}}
f(x)\geq4$, $\sigma(g_{X})\leq\sigma(f\upharpoonleft G_{X})+\frac{1}{4}%
{\textstyle\sum_{x\in X}}
f(x)$.

\item Suppose $X\in\mathcal{A}_{4}-\mathcal{A}_{4}^{\prime}$. Then $G_{X}$ is
a subgraph of some graph $G(k_{i})$, which will be considered next, when we
consider $\mathcal{E}$.

\item Each $X\in\mathcal{E}$ belongs to some $G(k_{i})$. Let $c_{k_{i}}$ be a
central vertex of $G(k_{i})$ and define the broadcast $g_{k_{i}}$ on
$G(k_{i})$ by%
\[
g_{k_{i}}(a)=\left\{
\begin{tabular}
[c]{ll}%
$1+%
{\textstyle\sum_{x\in V(G(k_{i}))}}
f(x)$ & if $a=c_{k_{i}}$\\
$0$ & otherwise.
\end{tabular}
\ \ \ \right.
\]
By (\ref{Eq_Gki}), $g_{k_{i}}$ dominates $G(k_{i})$ and, since $%
{\textstyle\sum_{x\in V(G(k_{i}))}}
f(x)\geq4$, $\sigma(g_{k_{i}})\leq\sigma(f\upharpoonleft G(k_{i}))+\frac{1}{4}%
{\textstyle\sum_{x\in V(G(k_{i}))}}
f(x)$.

\item Finally, let $X\in(\mathcal{A}_{1}\cup\mathcal{A}_{2}\cup\mathcal{A}%
_{3})-\mathcal{E}$. In each case we define a dominating broadcast $g_{X}$ on
$G_{X}$ such that $\sigma(g_{X})=\sigma(f\upharpoonleft G_{X})$.

\begin{enumerate}
\item If $X=\{x\}\in\mathcal{A}_{1}$, let $g_{X}(x)=f(x)=1$ and $g_{X}(a)=0$
for $a\in V(G_{X})-\{x\}$. Then $g_{X}$ dominates $G_{X}=G[N[x]]$.

\item Suppose $X=\{x\}\in\mathcal{A}_{2}\cup\mathcal{A}_{3}$. Each $u\in
U_{f}$ that annihilates $x$ is adjacent to all vertices in $\operatorname{PB}%
_{f}(x)=B_{f}(x)$, which induces a complete subgraph of $G$ (since
$X\notin\mathcal{D}$). In this case let $c_{X}$ be a vertex adjacent to $x$ on
an $x-u$ geodesic and define $g_{X}$ by $g_{X}(c_{X})=f(x)$ and $g_{X}(a)=0$
otherwise. If no $u\in U_{f}$ annihilates $x$, define $g_{X}$ by
$g_{X}(x)=f(x)$ and $g_{X}(a)=0$ otherwise. In either case $g_{X}$ dominates
$G_{X}$.

\item If $X=\{x,y\}\in\mathcal{A}_{2}$, or $X=\{x,y\}\in\mathcal{A}_{3}$ and
$d(x,y)=2$, or $X=\{x,y\}\in\mathcal{A}_{3}$, $d(x,y)=3$ and no vertex in
$U_{f}$ annihilates $y$, let $c_{X}$ be a central vertex of $G_{X}$ and define
$g_{X}$ by $g_{X}(c_{X})=f(x)+f(y)$ and $g_{X}(a)=0$ otherwise. Then $g_{X}$
dominates $G_{X}$.

\item Suppose $X=\{x,y\}\in\mathcal{A}_{3}$ and $d(x,y)=3$, where $f(x)=1$ and
$f(y)=2$, and some vertex in $U_{f}$ annihilates $y$. Since $X\notin%
\mathcal{C}$, some $y^{\prime}\in\operatorname{PB}_{f}(y)$ is adjacent to a
vertex $x^{\prime}\in B_{f}(x)$. Define $g_{X}$ by $g_{X}(y^{\prime})=3$, and
$g_{X}(a)=0$ otherwise.

\item Suppose $X=\{x,y,z\}\in\mathcal{A}_{3}$. In each case $c_{X}=y$ is a
central vertex of $G_{X}$; define $g_{X}$ by $g_{X}(c_{X})=3$, and
$g_{X}(a)=0$ otherwise. \label{defg}
\end{enumerate}
\end{enumerate}

Let $g$ be the union of all the functions $g_{X}$ and $g_{k_{i}}$ (seen as
sets of ordered pairs) defined above. Since the subgraphs considered above are
disjoint, $g$ is well defined. Also, since each $u\in U_{f}$ annihilates some
$x\in V_{f}^{+}$, each such $u$ belongs to some graph $G_{X}$ or $G(k_{i})$,
hence $g$ dominates $G$. Finally, as indicated in (1) to (4) above,
$\sigma(g)\leq\frac{5}{4}\sigma(f)$. In particular, choosing $f$ to be an
$\operatorname{ir}_{b}$ broadcast of $G$, it follows that $\gamma_{b}%
(G)\leq\frac{5}{4}\operatorname{ir}_{b}(G)$.~$\blacksquare$

\bigskip%
\begin{figure}[pb]%
\centering
\includegraphics[
height=1.8723in,
width=4.4797in
]%
{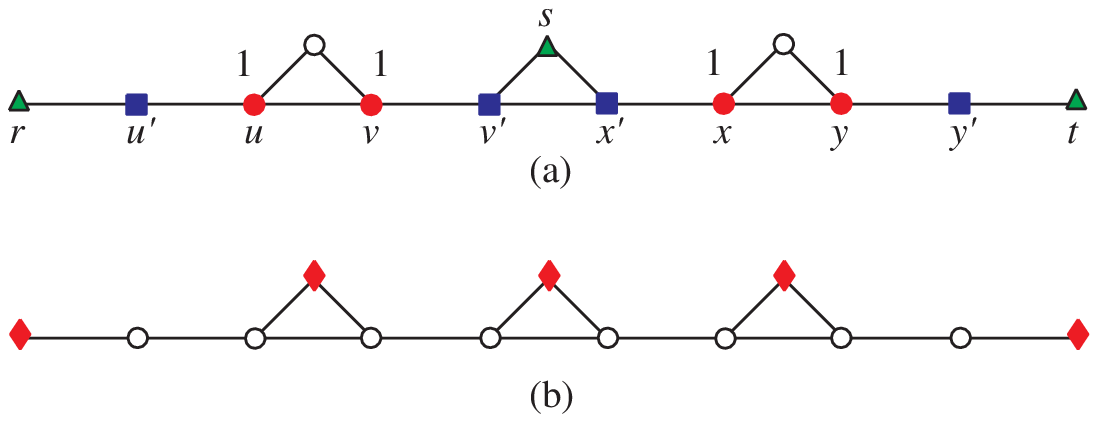}%
\caption{A graph $G$ with $\operatorname{ir}_{b}(G)=4$ and $\gamma
_{b}(G)=\operatorname{mp}(G)=5$}%
\label{Fig_ir_less_gamma}%
\end{figure}

As an example of a graph that attains equality in the bound in Theorem
\ref{Thm_irb-gammab_bound}, consider the graph $G$ in Figure
\ref{Fig_ir_less_gamma}, in which $f(w)=1$ if $w\in\{u,v,x,y\}$ (red solid
vertices in Figure \ref{Fig_ir_less_gamma}(a)) and 0 otherwise. This graph
appears in \cite{Ahmadi} as an example of a graph with $\operatorname{ir}%
_{b}<\gamma_{b}$. The unique vertex in the private boundary of each vertex
$w\in V_{f}^{+}$ is the blue square vertex $w^{\prime}$, and the three green
triangular vertices $r,\ s$ and $t$ are undominated. Clearly, Theorem
\ref{ThmMax_irr}$(i)$ holds for $V-V_{f}^{+}$. It is also easy to verify that
Theorem \ref{ThmMax_irr}$(ii)$ holds for each vertex in $V_{f}^{+}=V_{f}%
^{\ast}$. Hence $\operatorname{ir}_{b}(G)\leq4$. But $\gamma_{b}(G)=5$:
broadcast with a strength of 5 from $s$; all vertices hear this broadcast. The
red diamond-shaped vertices in Figure \ref{Fig_ir_less_gamma}(b) form a
multipacking of cardinality 5, thus confirming that $\gamma_{b}(G)=5$. By
Theorem \ref{Thm_irb-gammab_bound}, $\operatorname{ir}_{b}(G)=4$.

\begin{figure}[ptb]%
\centering
\includegraphics[
height=1.3007in,
width=3.9885in
]%
{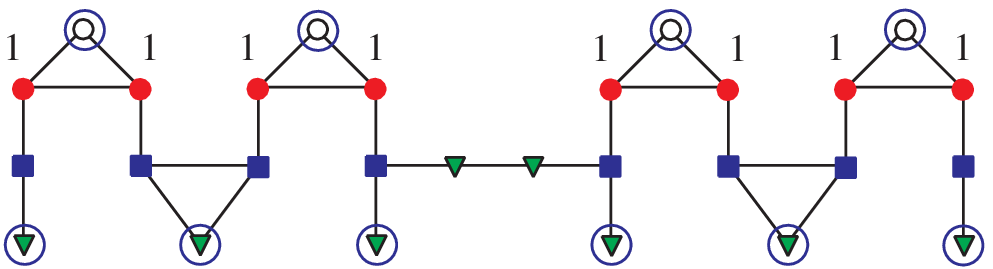}%
\caption{The graph $G_{2}$ with $\gamma_{b}(G_{2})=10$ and $\operatorname{ir}%
_{b}(G_{2})=8$}%
\label{Fig_ir_less_gamma2}%
\end{figure}

That the difference $\gamma_{b}-\operatorname{ir}_{b}$ can be arbitrarily
large for connected graphs can be seen by joining several copies of the graph
$G$ in Figure \ref{Fig_ir_less_gamma} linearly as shown in Figure
\ref{Fig_ir_less_gamma2}. If $G_{k}$ is the graph obtained by joining $k$
copies of $G$, then $\gamma_{b}(G_{k})=5k$ (the circled vertices form a
multipacking of cardinality $5k$) and $\operatorname{ir}_{b}(G_{k})=4k$ (by
Theorems \ref{Thm_irb-gammab_bound} and \ref{ThmMax_irr}). The graphs $G_{k}$
thus form an infinite class of graphs for which $\gamma_{b}/\operatorname{ir}%
_{b}=\frac{5}{4}$. (A larger class of graphs with this property can be
obtained by replacing the square vertices and circled vertices with cliques of
arbitrary order.) 

\begin{figure}[ptb]%
\centering
\includegraphics[
height=1.8239in,
width=6.4532in
]%
{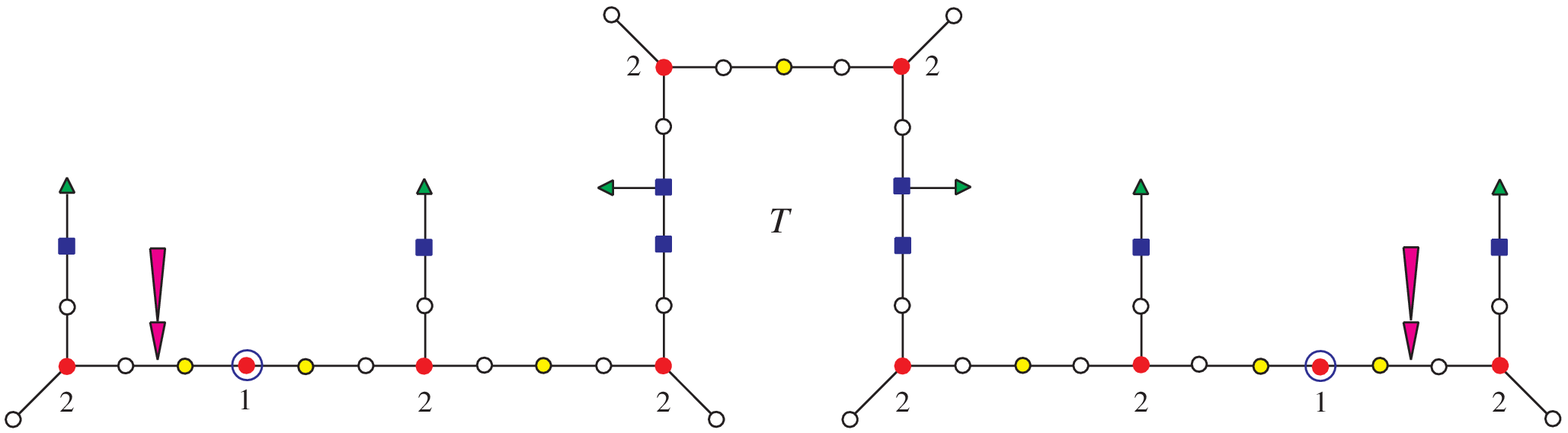}%
\caption{A tree $T$ with $\operatorname{ir}_{b}(T)=18$ and $\gamma_{b}(T)=19$}%
\label{Fig_irb-tree}%
\end{figure}

Consider the tree $T$ displayed in Figure \ref{Fig_irb-tree}. A maximal
irredundant broadcast $f$ with $\sigma(f)=18$ is shown, hence
$\operatorname{ir}_{b}(T)\leq18$. It can be shown that $\operatorname{ir}%
_{b}(T)=18$. Using the formula in \cite{HM} one easily obtains that
$\gamma_{b}(T)=19$. We believe $T$ to be the smallest tree for which
$\operatorname{ir}_{b}<\gamma_{b}$. If $T_{k}$ is the tree constructed from
$T$ similar to the construction of $G_{k}$ in Figure \ref{Fig_ir_less_gamma2}
from $G$, then $\operatorname{ir}_{b}(T_{k})\leq18k$ and $\gamma_{b}%
(T_{k})=19k$, hence $\gamma_{b}-\operatorname{ir}_{b}$ can also be arbitrarily
large for trees.%
%

\section{Upper Broadcast Domination and Irredundance}

\label{SecUpper}In this section we briefly discuss the relationships between
the upper domination, broadcast domination, irredundance and broadcast
irredundance numbers.%
\begin{figure}[pb]%
\centering
\includegraphics[
height=1.9147in,
width=1.7988in
]%
{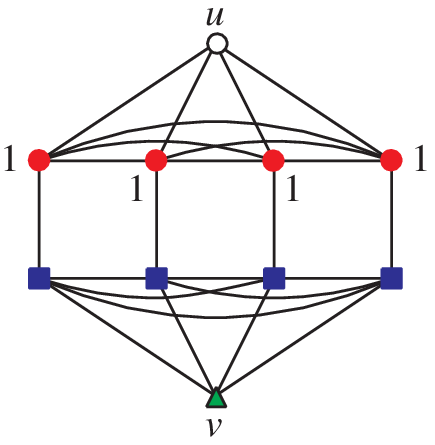}%
\caption{The graph $G_{4}$ with $\Gamma(G_{4})=2<\Gamma_{b}(G_{4}%
)=3<\operatorname{IR}(G_{4})=\operatorname{IR}_{b}(G_{4})=4$}%
\label{FigIR_Gamma}%
\end{figure}

The two rightmost inequalities in Corollary \ref{Cor_gammab-irb} can be
strict. For $r\geq3$, let $G_{r}$ be the graph obtained by joining two copies
of $K_{r+1}$ by $r$ independent edges. The graph $G_{4}$ is illustrated in
Figure \ref{FigIR_Gamma}. The red vertices (solid circles) form a maximal
irredundant set $S$ as well as a maximal irredundant broadcast $f$ with
$V_{f}^{+}=V_{f}^{1}=S$. Broadcasting to both $u$ and $v$ either requires a
broadcast from $u$ or $v$ with cost $3$, or a broadcast from a solid circle or
square with cost $2$, or a broadcast from one vertex in each $K_{5}$ with a
strength of 1. Hence $\Gamma(G_{4})=2<\Gamma_{b}(G_{4})=3<\operatorname{IR}%
(G_{4})=\operatorname{IR}_{b}(G_{4})=4$, and, in general, $\Gamma
(G_{r})=2<\Gamma_{b}(G_{r})=3<\operatorname{IR}(G_{r})=\operatorname{IR}%
_{b}(G_{r})=r$. Hence we see that $\operatorname{IR}_{b}/\Gamma_{b}$ is
unbounded. 

The difference $\operatorname{IR}_{b}-\Gamma_{b}$ can also be arbitrary for
trees, which is different from the situation for $\Gamma$ and
$\operatorname{IR}$: Cockayne, Favaron, Payan and Thomason \cite{CFPT} showed
that if $G$ is bipartite, then $\alpha(G)=\Gamma(G)=\operatorname{IR}(G)$,
where $\alpha(G)$ denotes the independence number of $G$. Figure
\ref{Fig_Trees}(a) shows a tree $T$ with $\operatorname{IR}_{b}(T)\geq14$; it
is not hard to show that equality holds. It follows from an argument in
\cite{Dunbar} for a supertree of $T$ that $\Gamma_{b}(T)=13$; see Figure
\ref{Fig_Trees}(b). Let $T=T_{1}$ and, for $k\geq2$, let $T_{k}$ be the tree
obtained from $k$ copies $S_{1},...,S_{k}$ of $T$ by joining the vertex
$w_{i}$ of $S_{i}$ to $u_{i+1}$ of $S_{i+1}$, $i=1,...,k-1$, as illustrated in
Figure \ref{Fig_Trees}. Then $\operatorname{IR}_{b}(T_{k})\geq14k$, while a
long and tedious argument, omitted here, shows that $\Gamma_{b}(T_{k})=13k$.
\emph{%
\begin{figure}[ptb]%
\centering
\includegraphics[
height=1.8723in,
width=6.2569in
]%
{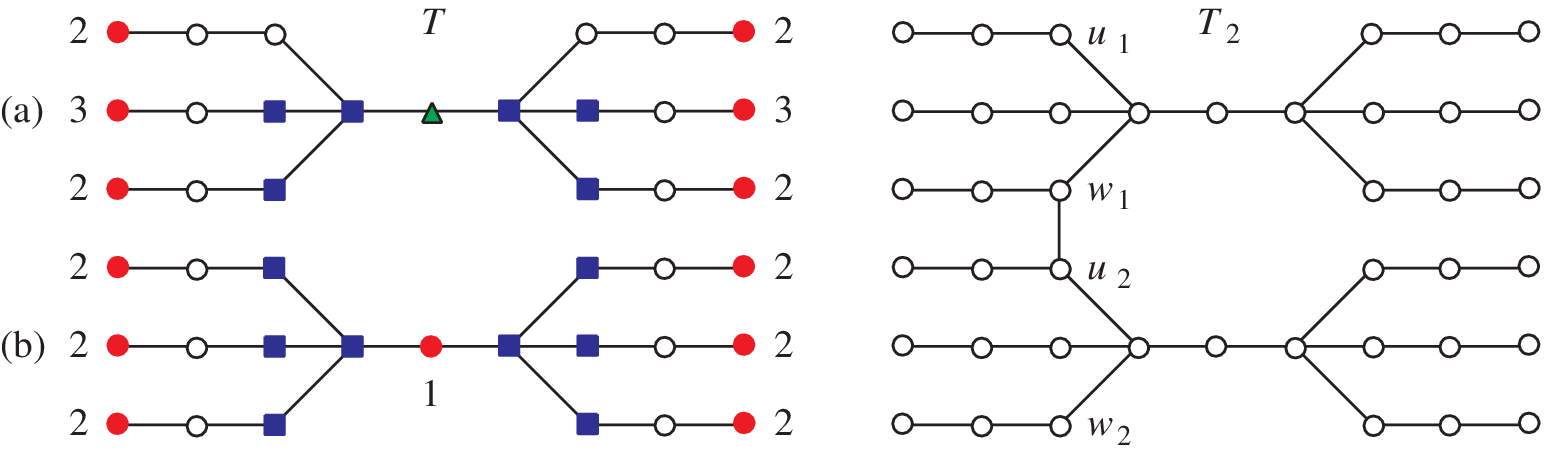}%
\caption{$\operatorname{IR}_{b}(T)=14$ and $\Gamma_{b}(T)=13$, while
$\operatorname{IR}_{b}(T_{2})\geq28$ and $\Gamma_{b}(T_{2})=26$}%
\label{Fig_Trees}%
\end{figure}
}

There also exist graphs such that $\Gamma_{b}$ and $\operatorname{IR}_{b}$ are
greater than $\operatorname{IR}$, as shown in Figure \ref{FigMax_ir}(c). The
grid $P_{m}\boksie P_{n}$ is bipartite, and it is well known that
$\alpha(P_{m}\boksie P_{n})=\left\lceil \frac{mn}{2}\right\rceil $ for all
$m,n$. Hence $\operatorname{IR}(P_{m}\boksie P_{n})=\left\lceil \frac{mn}%
{2}\right\rceil $. Figure \ref{FigMax_ir}(c) shows an irredundant and
dominating broadcast of $P_{4}\boksie P_{4}$ of cost $12>8=\alpha
(P_{4}\boksie
P_{4})$. Similarly, $\operatorname{IR}_{b}(P_{m}\boksie P_{n})\geq\Gamma
_{b}(P_{m}\boksie P_{n})\geq m(n-1)>\left\lceil \frac{mn}{2}\right\rceil
,\ m,n\geq3$.

\section{Open Problems}

\label{SecOpen}We conclude the paper with a list of open problems and two conjectures.

\begin{problem}
\label{Prob_Lower}Characterize the classes of graphs for which $\gamma
_{b}=\frac{5}{4}\operatorname{ir}_{b}$ and $\gamma_{b}=\operatorname{ir}_{b}$, respectively.
\end{problem}

We make the following two conjectures concerning Problem \ref{Prob_Lower}.

\begin{conjecture}
\label{Con5/4}For any graph $G$, $\gamma_{b}(G)=\frac{5}{4}\operatorname{ir}%
_{b}(G)$ if and only if each $\operatorname{ir}_{b}$-broadcast $f$ of $G$
satisfies $V_{f}^{+}=V_{f}^{1}$ and $G[V_{f}^{+}]=rK_{2}$.
\end{conjecture}

\begin{conjecture}
\label{Con=}If $G$ has an $\operatorname{ir}_{b}$-broadcast $f$ such that
$f(v)\geq2$ for all $v\in V_{f}^{+}$, then $\operatorname{ir}_{b}%
(G)=\gamma_{b}(G)$.
\end{conjecture}

The condition in Conjecture \ref{Con=} is not necessary for the two parameters
to be equal (consider complete graphs).

\begin{problem}
Determine the maximum ratio $\gamma_{b}/\operatorname{ir}_{b}$ for trees. The
tree in Figure $\ref{Fig_irb-tree}$ shows that $\max_{T\text{ a tree}}%
\{\gamma_{b}(T)/\operatorname{ir}_{b}(T)\}\geq\frac{19}{18}>1$.
\end{problem}

\begin{problem}
Compare $\operatorname{ir}$ and $\operatorname{ir}_{b}$. For example, is it
true that $\operatorname{ir}_{b}(T)\leq\operatorname{ir}(T)$ for all trees
$T$? Characterize graphs or trees for which $\operatorname{ir}_{b}%
=\operatorname{ir}$.
\end{problem}

\begin{problem}
Investigate $\Gamma_{b}$ and $\operatorname{IR}_{b}$. For example,\vspace
{-0.08in}

\begin{enumerate}
\item[$(i)$] characterize graphs for which $\Gamma_{b}=\operatorname{IR}_{b}%
$;\vspace{-0.08in}

\item[$(ii)$] determine the maximum ratio $\operatorname{IR}_{b}/\Gamma_{b}$
for bipartite graphs/trees. (The ratio is unbounded for general graphs, as
mentioned in Section \ref{SecUpper}.)
\end{enumerate}
\end{problem}

\begin{problem}
How do $\Gamma_{b}$ and $\operatorname{IR}_{b}$ compare to $\Gamma$ and
$\operatorname{IR}$ for various graph classes? What are the maximum ratios
$\Gamma_{b}/\Gamma$ and $\operatorname{IR}_{b}/\operatorname{IR}$?
\end{problem}

\begin{problem}
Construct algorithms to determine $\operatorname{ir}_{b}$, $\operatorname{IR}%
_{b}$, $\Gamma_{b}$. What is the complexity of determining $\operatorname{ir}%
_{b}$, $\operatorname{IR}_{b}$, $\Gamma_{b}$? (Heggernes and Lokshtanov
\emph{\cite{HL}} showed that determining $\gamma_{b}$ is $\mathcal{O}(n^{8})$.)
\end{problem}

\begin{problem}
As mentioned in Section $\ref{Sec_Dom}$, $\Gamma_{b}(G)\geq\operatorname{diam}%
(G)$ for all graphs $G$. Call a graph $G$ a \emph{diametrical graph} if
$\Gamma_{b}(G)=\operatorname{diam}(G)$. Characterize diametrical
trees.\label{prob}
\end{problem}

\begin{problem}
Study criticality concepts for $\operatorname{ir}_{b}$, $\operatorname{IR}%
_{b}$, $\gamma_{b}$, $\Gamma_{b}$.
\end{problem}

\label{refs}

\bigskip

\begin{thebibliography}{99}                                                                                               %


\bibitem {Ahmadi}D.~Ahmadi, G.~H.~Fricke, C.~Schroeder, S.~T.~Hedetniemi,
R.~C.~Laskar, Broadcast irredundance in graphs. \emph{Congr.~Numer.}%
~\textbf{224} (2015), 17--31.

\bibitem {BC}B.~Bollob\'{a}s, E.~J.~Cockayne, Graph-theoretic parameters
concerning domination, independence, and irredundance, \emph{J.~Graph Theory}
\textbf{3} (1979), 241--249.

\bibitem {Bouch}I.~Bouchemakh, M.~Zemir, On the broadcast independence number
of grid graph, \emph{Graphs Combin.} \textbf{30} (2014), 83--100.

\bibitem {BMT}R.~C.~Brewster, C.~M.~Mynhardt, L.~E.~Teshima, New bounds for
the broadcast domination number of a graph, \emph{Cent. Eur. J. Math.}
\textbf{11}(7) (2013), 1334--1343.

\bibitem {CLZ}G.~Chartrand, L.~Lesniak, P.~Zhang, \emph{Graphs \&\ Digraphs},
Chapman and Hall/CRC, Boca Raton, 2016.

\bibitem {CHM}E.~J.~Cockayne, S.~T.~Hedetniemi, D.~J.~Miller, Properties of
hereditary hypergraphs and middle graphs, \emph{Canad.~Math.~Bull.}
\textbf{21} (1978), 461--468.

\bibitem {CFPT}E.~J.~Cockayne, O.~Favaron, C.~Payan, A.~G.~Thomason,
Contributions to the theory of domination, independence and irredundance in
graphs, \emph{Discrete Math.~}\textbf{33 }(1981), 249--258.

\bibitem {DDH}J.~Dabney, B.~C.~Dean, S.~T.~Hedetniemi, A linear-time algorithm
for broadcast domination in a tree, \emph{Networks. }\textbf{53} (2009), 160--169.

\bibitem {Dunbar}J.~Dunbar, D.~Erwin, T.~Haynes, S.~M.~Hedetniemi,
S.~T.~Hedetniemi, Broadcasts in graphs, \emph{Discrete Applied Math.}
\textbf{154} (2006), 59-75.

\bibitem {Ethesis}D.~Erwin, \emph{Cost domination in graphs}. Doctoral
dissertation, Western Michigan University, 2001.

\bibitem {HL}P.~Heggernes, D.~Lokshtanov, Optimal broadcast domination in
polynomial time. \emph{Discrete Math.} \textbf{36} (2006), 3267--3280.

\bibitem {Herke}S.~Herke, \emph{Dominating broadcasts in graphs}, Master's
thesis, University of Victoria, 2009. http://hdl.handle.net/1828/1479

\bibitem {HM}S.~Herke, C.~M.~Mynhardt, Radial Trees. \emph{Discrete Math.
}\textbf{309} (2009), 5950--5962.

\bibitem {LauraT}L. Teshima, \emph{Broadcasts and multipackings in graphs},
Master's thesis, University of Victoria, 2012. http://hdl.handle.net/1828/4341

\bibitem {MT}C.~M.~Mynhardt, L.~E.~Teshima, Broadcasts and multipackings in
trees, \emph{Utilitas Math.}, to appear.
\end{thebibliography}
\end{document}